\documentclass[reqno,a4paper,12pt]{amsart}
\usepackage[top=32mm, bottom=32mm, left=32mm, right=32mm]{geometry}
\usepackage{mathptmx}
\usepackage{amssymb,latexsym}
\usepackage{enumerate}
\usepackage{color}
\usepackage[mathcal]{euscript}
\usepackage{mathrsfs}
\usepackage{graphicx}
\usepackage[colorlinks=true,citecolor=blue]{hyperref}

\theoremstyle{plain}
\newtheorem{theorem}{Theorem}[section]
\newtheorem*{LDT-thm}{Large deviations theorem}
\newtheorem*{BET-thm}{Birkhoff Ergodic Theorem}
\newtheorem{lemma}{Lemma}[section]
\newtheorem{proposition}{Proposition}[section]
\newtheorem{corollary}{Corollary}[section]
\newtheorem{claim}{Claim}

\theoremstyle{definition}

\newtheorem{example}{Example}[section]
\theoremstyle{question}
\newtheorem{question}{Question}

\theoremstyle{remark}
\newtheorem{remark}{Remark}

\numberwithin{equation}{section}

\numberwithin{equation}{section}

\begin{document}

\title[]{On the large deviations theorem of weaker types}

\author[X. Wu]{Xinxing Wu}
\address[X. Wu]{School of Sciences, Southwest Petroleum University, Chengdu, Sichuan, 610500, People's
Republic of China}
\email{wuxinxing5201314@163.com}

\author[X. Wang]{Xiong Wang}
\address[X. Wu;X. Wang]{Institute for Advanced Study, Shenzhen University, Nanshan District Shenzhen, Guangdong, China}
\email{wangxiong8686@szu.edu.cn}
\author[G. Chen]{Guanrong Chen}
\address[G. Chen]{Department of Electronic
Engineering, City University of Hong Kong, Hong Kong SAR, People's
Republic of China}
\email{eegchen@cityu.edu.hk}
\thanks{Xinxing Wu was supported by
the Scientific Research Fund of Sichuan Provincial Education
Department (No. 14ZB0007) and the
scientific research starting project of SWPU}
\thanks{Xiong Wang (Corresponding Author) was supported by the National Natural Science Foundation of China (No.11547117) and Seed Funding from Scientific and Technical Innovation Council of Shenzhen Government.}
\thanks{Guanrong Chen was supported by the Hong Kong Research Grants Council under
GRF Grant CityU 11201414.}

\subjclass[2010]{Primary 54H20, 37B99; Secondary 54B20, 54H20, 37B20, 37B40.}
\date{\today}


\keywords{Large deviations theorem, ergodic, transitivity, sensitivity, equicontinuous.}
\maketitle
\begin{abstract}
  In this paper, we introduce
  the concepts of the large deviations theorem of weaker types, i.e.,
  type I, type I', type II, type II', type III,
  and type III', and present a systematic study of the ergodic and chaotic properties
  of dynamical systems satisfying the large deviations theorem of various types.
  Some characteristics of the ergodic measure
  are obtained and then applied to prove that every dynamical system
  satisfying the large deviations theorem of type I' is ergodic, which is
  equivalent to the large deviations theorem of type II' in this regard,
  and that every uniquely ergodic dynamical system restricted on its support
  satisfies the large deviations theorem.
  Moreover, we prove that every dynamical system satisfying the large deviations
  theorem of type III is an $E$-system. Finally, we show that
  a dynamical system satisfying the central limit theorem,
  introduced in [Y. Niu, Y. Wang, Statist. Probab. Lett.,
   {\bf 80} (2010), 1180--1184], does not exist. 
\end{abstract}


\section{Introduction}

A {\it dynamical system} is a pair $(X, T)$, where $X$ is a compact
metric space with a metric $d$ and $T: X\rightarrow X$ is a
continuous map. Sharkovsky's amazing discovery \cite{Sharkovskii1964}, as well
as Li and Yorke's famous work which introduced the concept of
`chaos' known as Li-Yorke chaos today
\cite{Li-Yorke1975}, have provoked the recent rapid advancement of research on
discrete chaos theory. At the same time, many
research works were devoted to the links between the topological and
stochastic properties of deterministic dynamical systems.

The topological approach on chaoticity tries to describe the
topological structure of a `chaotic region' \cite{Akin-Kolyada2003,
Block-Coppel1992, Devaney1989, Huang-Ye2001, Huang-Ye2002,
Kolyada-Snoha1997, Li-Yorke1975, Moothathu2007}. The essence of Li-Yorke chaos is the
existence of uncountable scrambled sets. Another well-known definition
of chaos was given by Devaney \cite{Devaney1989}, according to which
a continuous map $T$ is said to be {\it chaotic in the sense of
Devaney}, if it satisfies the following three properties: With notation $\mathbb{N}=\{1, 2, \ldots\}$,
$\mathbb{Z}^{+}=\{0, 1, 2, \ldots\}$,
\begin{enumerate}
\item $T$ is {\it topologically transitive}, i.e., for every pair of
nonempty open sets $U, V\subset X$, there exists
$n\in \mathbb{Z}^{+}$ such that
$T^{n}(U)\cap V\neq \O$;

\item The set of periodic points of $T$ is dense in $X$;

\item\label{D-3} $T$ has {\it sensitive dependence on initial conditions}
(briefly, {\it is sensitive}), i.e., there exists $\varepsilon >0$
such that for any $x\in X$ and any neighborhood $U$ of $x$,
there exist $y\in U$ and $n\in \mathbb{Z}^{+}$ satisfying
$d(T^{n}(x), T^{n}(y))>\varepsilon$.
\end{enumerate}

Banks et al. \cite{Banks1992} proved that
every topologically transitive map whose periodic points are dense
in $X$ has sensitive dependence on initial conditions, which implies that
the above condition (\ref{D-3}) is redundant, while Huang and Ye \cite{Huang-Ye2002} showed that every
topologically transitive map containing a periodic point is chaotic
in the sense of Li-Yorke. Most importantly, sensitive
dependence on initial conditions is widely understood as a key
ingredient of chaos and was popularized by the meteorologist Lorenz,
which is commonly known as the so-called `butterfly effect'.

The stochastic approach is devoted to characterizing the dynamics of
a deterministic dynamical system through its stochastic properties
by using tools from ergodic theory, functional analysis and
spectral theory \cite{Baladi2000,Liverani1995a,Liverani1995b,O1952,
Viana1997,Walters1972,Walters1982}. Meanwhile, the links between the two
approaches are being gradually studied
\cite{Abraham2002,Gu2007, He2004,Lardjane2006,Li2013,
Niu-Wang2010,Wu2005,Xu2004,Zhao2002}. For example, Abraham et al.
\cite{Abraham2002} gave some sufficient conditions (in terms of topology
and ergodicity) on a measure-preserving dynamical system defined on a
nontrivial metric space $(X, \mathcal{B}(X), \mu)$ endowed with a
Borel probability measure, to ensure the sensitivity property or
cofinite sensitivity property. In 2004, He et al. \cite{He2004}
proved that for a measure-preserving  transformation $T$ on $(X, \mathcal{B}(X),
\mu)$, if $\mathrm{supp}(\mu)=X$ and $T$ is weakly mixing, then $T$ has
sensitive dependence on initial conditions and this also holds
for measure-preserving semiflows. In 2006, Lardjane
\cite{Lardjane2006} complemented the main results in \cite{
Abraham2002, Wu2005, Xu2004} on the links between several
topological and stochastic properties of dynamical systems and
proved that if $\mathrm{supp}(\mu)=X$ and $T$ is mixing
(not-necessarily measure-preserving), then $T$ has sensitive
dependence on initial conditions.

It is commonly known that the large deviations theorem and the
central limit theorem coming from probability theory are two
of the most remarkable results in all fields of mathematics
especially in probability and
statistics. They were successfully applied to dynamical
systems \cite{Lardjane2006,Li2013,Niu2009-2,Niu-Wang2010,Wu2005,WC2014,WC2015,Xu2004}.
The former describes the oscillation of the time average
$(1/ n)\sum_{i=0}^{n-1}\varphi\circ T^{i}(x)$ around the spatial
average $\int_{X}\varphi\mathrm{d}\mu$ and the latter describes the
rate of its convergence. In 2007, Gu \cite{Gu2007} extended the
results obtained by Wu et al. \cite{Wu2005} and showed that
 a dynamical system satisfying the large deviations theorem is
topologically ergodic. Moreover, if it is strongly topologically
ergodic, then it has sensitive dependence on initial conditions
(see \cite[Theorem 3.1, Theorem 4.1]{Gu2007}).
Lately, Niu \cite{Niu2009-2} proved that a dynamical system satisfying
the large deviations theorem has an equicontinuous point if and
only if it is both minimal and equicontinuous.
Then, Li \cite{Li2013} introduced the concept of ergodic
sensitivity, which is a stronger form of sensitivity, and showed
that a strongly topologically ergodic system satisfying
the large deviations theorem is ergodically sensitive. More recently, we
\cite{WC2015} proved that a dynamical system satisfying the
large deviations theorem is ergodic.

Based on the results in \cite{Gu2007,Lardjane2006,Li2013,Niu2009-2,Niu-Wang2010,WC2015},
our objective here is to use the methods
of ergodic theory and topological dynamics to further investigate the relations
between the large deviations theorem, the ergodic properties and
the chaotic behaviors of dynamical systems.


\section{Preliminaries}

\subsection{Topological dynamics}

For $U, V\subset X$, define the {\it return time set from $U$
to $V$} as $ N(U, V)=\{n\in \mathbb{Z}^{+}:
T^{n}(U)\cap V\neq \O \}$. In particular,
$N(x, V)=\left\{n\in \mathbb{Z}^{+}:
T^{n}(x)\in V\right\}$ for $x\in X$.

Let $\mathcal{P}$ be the
collection of all subsets of $\mathbb{Z}^{+}$. A collection
$\mathscr{F}\subset\mathcal{P}$ is called a {\it Furstenberg family}
if it is hereditary upwards, i.e., $F_{1}\subset F_{2}$ and
$F_{1}\in\mathscr{F}$ imply $F_{2}\in\mathscr{F}$. A family $\mathscr{F}$
is {\it proper} if it is a proper subset of $\mathcal{P}$, i.e.,
neither empty nor the whole $\mathcal{P}$. It is easy to see
that $\mathscr{F}$ is proper if and only if $\mathbb{Z}^{+}\in \mathscr{F}$
and $\O\notin \mathscr{F}$. All the families considered below are assumed to be
proper. For a Furstenberg family $\mathscr{F}$, denote $\Delta(\mathscr{F})=\{F-F: F\in \mathscr{F}\}$,
where $F-F=\{i-j\in \mathbb{Z}^{+}: i, j\in F\}$.

For $A\subset \mathbb{Z}^{+}$, define
$$
\overline{d}(A)=\limsup_{n\rightarrow +\infty}\frac{1}{n}\left| A\cap
[0, n-1]\right|
\text{ and }
\underline{d}(A)=\liminf_{n\rightarrow +\infty}\frac{1}{n}\left| A\cap
[0, n-1]\right|.
$$
Then, $\overline{d}(A)$ and $\underline{d}(A)$ are {\it the upper
density} and {\it the lower density} of $A$, respectively.

Similarly, define the \textit{upper Banach density}  and the \textit{lower Banach density}
of $A$ as
$$
\mathrm{BD}^*(A)=\limsup_{|I| \rightarrow +\infty} \frac{|A \cap I|}{|I|} \text{ and }
\mathrm{BD}_*(A)=\liminf_{|I| \rightarrow +\infty} \frac{|A \cap I|}{|I|},
$$
where $I$ is over all non-empty finite intervals of
$\mathbb{Z}^{+}$.

A subset $S$ of $\mathbb{Z}_+$ is {\it syndetic} if it has a bounded
gap, i.e., if there is $N \in \mathbb{N}$ such that $\{i, i+1, \ldots,
i+N\} \cap S \neq \emptyset$ for every $i \in \mathbb{Z}^+$; $S$ is
{\it thick} if it contains arbitrarily long runs of positive
integers, i.e., for every $n \in \mathbb{N}$ there exists some $a_n
\in \mathbb{Z}^+$ such that $\{a_n, a_n+1, \ldots, a_n+n\} \subset
S$. The set of all thick subsets of $\mathbb{Z}^{+}$, all syndetic
subsets of $\mathbb{Z}^{+}$, all subsets of $\mathbb{Z}^{+}$ with positive
upper density, all subsets of $\mathbb{Z}^{+}$ with upper density
equal to 1, and all subsets of $\mathbb{Z}^{+}$ with positive
upper Banach density, are denoted by $\mathscr{F}_{t}$, $\mathscr{F}_{s}$,
$\mathscr{F}_{pud}$, $\mathscr{F}_{ud1}$, and $\mathscr{F}_{pubd}$, respectively. Clearly,
all of them are Furstenberg families.

For a Furstenberg family $\mathscr{F}$, a dynamical system
is called {\it $\mathscr{F}$-transitive} if $N(U, V)\in \mathscr{F}$
for every pair of nonempty open subsets $U, V\subset X$.
The $\mathscr{F}_{pud}$-transitivity and $\mathscr{F}_{ud1}$-transitivity
are called topological ergodicity and strongly topological
ergodicity respectively in \cite{Gu2007, Li2013}. Clearly, $\mathscr{F}_{s}$-transitivity
is stronger that $\mathscr{F}_{pud}$-transitivity.

A dynamical system $(X, T)$ is {\it totally transitive} if $(X, T^n)$ is transitive
for each $n \in \mathbb{N}$; and it is {\it (topologically) weakly mixing} if $(X \times X, T \times
T)$ is transitive. It is well known that $(X, T)$ is weakly mixing if and only if
it is $\mathscr{F}_{t}$-transitive (see \cite{Fur-book,Fur1982}). An $x\in X$ is
\textit{a transitive point} if its orbital closure
$\overline{\textrm{orb}(x, T)}=X$. Let $\textrm{Trans}(T)$ be the set of transitive points.
Then, the orbit closure of a recurrent point is transitive.

A dynamical system $(X, T)$ is {\it minimal} if every orbit
under $T$ is dense in $X$. It is easy to see that $(X, T)$ is a minimal
system if and only if $X$ has no proper, nonempty, closed invariant
subset. A point $x\in X$ is called an {\it equicontinuity point} of $T$
if, for any $\varepsilon>0$, there exists $\delta>0$ such that for
any $y\in X$ with $d(x, y)< \delta$ and any $n\in \mathbb{Z}^{+}$, one has
$d(T^{n}(x), T^{n}(y))<\varepsilon$. A dynamical system
$(X, T)$ is {\it equicontinuous} if every $x\in X$ is an
equicontinuous point of $T$; $(X, T)$ is {\it almost equicontinuous} if it
is a transitive dynamical system admitting an equicontinuity point. By compactness,
it can be verified that $(X, T$) is equicontinuous
when the sequence $\{T^n: n\in \mathbb{Z}^{+}\}$
is uniformly equicontinuous.

For $U \subset X$ and $\varepsilon>0$, let
$$
N(U, \varepsilon)=\{n\in \mathbb{Z}^{+}: \mathrm{diam}T^{n}(U)>\varepsilon\}.
$$

It is easy to see that a dynamical system $(X, T)$ is sensitive if and only
if there exists $\varepsilon>0$ such that, for any nonempty open subset
$U\subset X$, $N(U, \varepsilon)\neq \O$.
For a dynamical system, Moothathu~\cite{Moothathu2007} initiated a
preliminary study of stronger forms of sensitivity formulated in
terms of some subsets of $\mathbb{Z}^{+}$, namely the
syndetical sensitivity and cofinite sensitivity. Recently, Li~\cite{Li2013}
introduced the concept of ergodic sensitivity. According to
Moothathu \cite{Moothathu2007} and Li \cite{Li2013},
a dynamical system $(X, T)$ is said to be
\begin{enumerate}
\item {\it ergodically sensitive} if there exists $\varepsilon>0$
  such that for any nonempty open subset $U\subset X$,
  $\overline{d}(N(U, \varepsilon))>0$;
\item {\it syndetically sensitive} if there exists $\varepsilon>0$
  such that for any nonempty open subset $U\subset X$,
  $N(U, \varepsilon)$ is syndetic;
\item {\it cofinitely sensitive} if there exists $\varepsilon>0$
  such that for any nonempty open subset $U\subset X$,
  $N(U, \varepsilon)$ is cofinite.
\end{enumerate}
Clearly, cofinite sensitivity is stronger than syndetical sensitivity,
which implies ergodic sensitivity. More results on sensitivity can be found in
\cite{HKZ,WWC2015,WOC}.

\subsection{Probability measure}
Let $\mathcal{B}(X)$ be the $\sigma$-algebra of Borel subsets of $X$,
$M(X)$ the set of Borel probability measures on $(X, \mathcal{B}(X))$,
and $M(X, T)$ the $T$-invariant ones. It is well known that $M(X)$ is a compact
metrisable space in the weak*-topology, and $M(X, T)$ is a nonempty closed
subset of $M(X)$. A measure-preserving transformation $T$ of
a probability space $(X, \mathcal{B}(X), \mu)$ is called
{\it ergodic} if the only members $B$ of $\mathcal{B}(X)$ with $T^{-1}(B)=B$ satisfy
$\mu(B)=0$ or $\mu(B)=1$. A probability measure $\mu\in M(X, T)$ is
called {\it ergodic} if the measure-preserving transformation $T$ of
the measure space $(X, \mathcal{B}(X), \mu)$ is ergodic. Let $E(X, T)$ be the set of all
ergodic measures in $M(X, T)$. If $M(X, T)$ consists of a single
point, then $(X, T)$ is said to be {\it uniquely ergodic}.

For any $\mu\in M(X)$, the set
$\{x\in X: \mu (U)>0 \text{ for any neighborhood }
U \text{ of } x\}$ is called the {\it support} of $\mu$, denoted by $\mathrm{supp}(\mu)$.
As every Borel probability measure $\mu$ is regular (see \cite[Theorem 6.1]{Walters1982}), it is easy to see
that $\mu(\mathrm{supp}(\mu))=1$ and $\mathrm{supp}(\mu|_{\mathcal{B}(\mathrm{supp}(\mu))})=\mathrm{supp}(\mu)$.

Let $C(X)$ denote the Banach space of continuous complex-valued functions on $X$ with the
supremum norm $\|\cdot\|$ and call each element of $C(X)$ an {\it observable}.
\begin{lemma} \label{lem:2-1}\cite[pp. 149]{Walters1982}
The following statements are equivalent:
\begin{enumerate}
  \item $\mu_n \rightarrow \mu$ in the weak*-topology.
  \item For each $\varphi\in C(X)$, $\int_{X}\varphi\mathrm{d}\mu_{n} \rightarrow \int_{X}\varphi \mathrm{d}\mu$ as $n\rightarrow +\infty$.
  \item For each closed subset $F\subset X$, $\limsup_{n\rightarrow+\infty}\mu_{n}(F)\leq \mu(F)$.
  \item For each open subset $U\subset X$, $\liminf_{n\rightarrow+\infty}\mu_{n}(U)\geq \mu(U)$.
\end{enumerate}
\end{lemma}

For $x \in X$, let $\delta_x \in M(X)$ denote the {\it Dirac point measure} of $x$, defined by
$$
\delta_x(A)=\left\{ \begin{array}{lc}
1, & x \in A, \\
0, & x \notin A.
\end{array}
\right.
$$

For the ergodic measure, the following result is well known. 
\begin{lemma}
\label{ergodic-support}
Let $(X, T)$ be a dynamical system and $\mu\in M(X, T)$. Then,
\begin{enumerate}
  \item $\mathrm{supp}(\mu)$ is a nonempty, closed, invariant subset of $X$.
  \item If $\mu$ is ergodic, then $(\mathrm{supp}(\mu), \mathcal{B}(\mathrm{supp}(\mu)),
  T, \mu)$ is ergodic and transitive.
  \item If $T$ is uniquely ergodic, then $(\mathrm{supp}(\mu), T)$ is uniquely ergodic
  and minimal.
\end{enumerate}
\end{lemma}

A dynamical system $(X, T)$ is called an {\it $E$-system} if it is topologically transitive
and there exists an invariant measure $\mu$ with a full support, i.e.,
$\mathrm{supp}(\mu)= X$. It is well known that every minimal system is an
$E$-system, and every $E$-system is $\mathscr{F}_{s}$-transitive (see \cite[Theorem 4.4]{GW2000}).

The main concern in the stochastic analysis
of a deterministic dynamical system is how to describe the
oscillations of the finite-time average
$$
\frac{1}{n}\sum_{i=0}^{n-1}\varphi \circ T^{i}(x)
$$
around their expected value $\int_{X}\varphi \mathrm{d}\mu$, where
$\varphi$ is an {\it observable},
i.e., $\varphi\in C(X)$. The Birkhoff Ergodic Theorem indicates that for an
ergodic measure-preserving transformation $T$ and for any $\varphi\in L^{1}(\mu)$,
$$
\lim_{n\rightarrow +\infty}\frac{1}{n}\sum_{i=0}^{n-1}\varphi \circ T^{i}(x)
=\int_{X}\varphi\mathrm{d}\mu, \ \mathrm{a.e.}
$$
\begin{BET-thm}
Suppose that $T$ is a measure-preserving transformation of a probability
space $(X, \mathcal{B}(X), \mu)$ and $\varphi\in L^{1}(\mu)$. Then,
$(1/ n)\sum_{i=0}^{n-1}\varphi \circ T^{i}(x)$ converges $\mathrm{a.e.}$ to a limit
function $\varphi^{*}\in L^{1}(\mu)$ such that $\varphi^{*}\circ T=\varphi^{*}$
and $\int_{X}\varphi^{*}\mathrm{d}\mu=\int_{X}\varphi\mathrm{d}\mu$.
\end{BET-thm}

When $T$ is uniquely ergodic, the following result shows that it admits much
stronger properties of these ergodic averages (see \cite[Theorem 6.19]{Walters1982}).
\begin{lemma}\cite[Theorem 6.19]{Walters1982}\label{Unique ergodic}
Let $(X, T)$ be a dynamical system. The following statements are equivalent:
\begin{enumerate}
  \item For every $\varphi\in C(X)$, $(1/ n)\sum_{i=0}^{n-1}\varphi\circ T^{i}(x)$ converges uniformly
  to a constant.
  \item For every $\varphi\in C(X)$, $(1/ n)\sum_{i=0}^{n-1}\varphi\circ T^{i}(x)$ converges pointwise
  to a constant.
  \item There exists $\mu\in M(X, T)$ such that, for all $\varphi\in C(X)$ and all $x\in X$,
  $$
  \lim_{n\rightarrow +\infty}\frac{1}{n}\sum_{i=0}^{n-1}\varphi\circ T^{i}(x)=\int_{X}\varphi\mathrm{d}\mu.
  $$
  \item $T$ is uniquely ergodic.
\end{enumerate}
\end{lemma}

Meanwhile, by the Birkhoff Ergodic Theorem, it can be verified that a measure-preserving
transformation $T$ of a probability space $(X, \mathcal{B}(X, \mu)$ is ergodic if and only
if for any $A, B\in \mathcal{B}(X)$,
$$
\lim_{n\rightarrow +\infty}\frac{1}{n}\sum_{i=0}^{n-1}\mu(A\cap T^{-i}(B))=\mu(A)\mu(B).
$$

\subsection{The large deviations theorem}
First, recall the original large deviations theorem in classic probability theory.
\begin{LDT-thm}\label{LDT-thm}
Let $X_{0}, X_{1}, \ldots$ be independent identically distributed random variables
taking values in $\mathbb{R}$, with average $\overline{X}=E(X_{n})<+\infty$, variance
$\sigma^{2}=E((X_{n}-\overline{X})^{2})\in (0, +\infty)$, and $E(e^{tX_{n}})\in (0, +\infty)$ for
every $t\in \mathbb{R}$. Then, for any $\varepsilon>0$, the probability $\mathbf{P}(n, \varepsilon)$
of
$$
\left|\frac{1}{n}\sum_{i=0}^{n-1}(X_{i}-\overline{X})\right|>\varepsilon
$$
converges to zero exponentially as $n\rightarrow +\infty$, in the sense that
$$
\limsup_{n\rightarrow +\infty}\frac{1}{n}\log \mathbf{P}(n, \varepsilon)<0.
$$
\end{LDT-thm}
Following the Birkhoff Ergodic Theorem and this large deviations theorem,
Wu et al. \cite{Wu2005} introduced the large deviations theorem for dynamical systems.

Let $(X, T)$ be a dynamical system and $\mu\in M(X)$.
An observable $\varphi$ is said to satisfy the {\it large deviations
theorem} for $(X, \mathcal{B}(X), T, \mu)$, or simply $(T, \mu)$,
if for any $\varepsilon>0$ there exists
$h(\varepsilon)>0$ such that
\begin{equation}\label{*}
\mu\left(\left\{x\in X: \left|
\frac{1}{n}\sum_{i=0}^{n-1}\varphi \circ T^{i}(x)
-\int_{X}\varphi\mathrm{d}\mu \right|>\varepsilon \right\}\right)\leq
e^{-nh(\varepsilon)}
\end{equation}
for all sufficiently large $n\in \mathbb{N}$. According to Wu et al. \cite{Wu2005},
a dynamical system $(X, \mathcal{B}(X), T, \mu)$ or $(T, \mu)$ is said to
satisfy the {\it large deviations theorem} 
if every observable
satisfies the large deviations theorem for $(T, \mu)$ and
$\mathrm{supp}(\mu)=X$.

To extend the large deviations theorem (LDT) modifying the convergence
in \eqref{*}, we now introduce some concepts
of large deviations theorem of weaker forms. We say that $(T, \mu)$ satisfies
\begin{enumerate}
  \item the {\it large deviations theorem of weak form} (WLDT) if, every
  observable satisfies the large deviations theorem.
  \item the {\it large deviations theorem of
  type I'} (LDT-I') if, for every observable $\varphi$ and any $\varepsilon>0$,
$$
\sum_{n=1}^{+\infty}\mu\left(\left\{x\in X: \left|
\frac{1}{n}\sum_{i=0}^{n-1}\varphi \circ T^{i}(x)
-\int_{X}\varphi\mathrm{d}\mu \right|>\varepsilon \right\}\right)< +\infty;
$$
  \item the {\it large deviations theorem of
  type II'} (LDT-II') if, for every observable $\varphi$ and any $\varepsilon>0$,
$$
\lim_{n \rightarrow +\infty}\mu\left(\left\{x\in X: \left|
\frac{1}{n}\sum_{i=0}^{n-1}\varphi \circ T^{i}(x)
-\int_{X}\varphi\mathrm{d}\mu \right|>\varepsilon \right\}\right)=0;
$$
  \item the {\it large deviations theorem of
  type III'} (LDT-III') if, for every observable $\varphi$ and any $\varepsilon>0$,
$$
\liminf_{n \rightarrow +\infty}\mu\left(\left\{x\in X: \left|
\frac{1}{n}\sum_{i=0}^{n-1}\varphi \circ T^{i}(x)
-\int_{X}\varphi\mathrm{d}\mu \right|>\varepsilon \right\}\right)=0.
$$
\end{enumerate}
If $(T, \mu)$ satisfies LDT-I' (resp., LDT-II', LDT-III')
and $\mathrm{supp}(\mu)=1$, then we say that $(T, \mu)$ satisfies the large deviations theorem
of type I (LDT-I) (resp., type II (LDT-II), type III (LDT-III)).

Clearly,
$$
\text{LDT} \Rightarrow \text{LDT-I} \Rightarrow \text{LDT-II} \Rightarrow
\text{LDT-III},
$$
and
$$
\text{WLDT} \Rightarrow \text{LDT-I'} \Rightarrow \text{LDT-II'}
\Rightarrow \text{LDT-III'}.
$$
Meanwhile, it can be verified that every trivial dynamical system satisfies LDT,
and that $(X, \mathcal{B}(X), T, \mu)$ satisfies LDT-I' (resp., LDT-II', LDT-III')
if and only if $(\mathrm{supp}(\mu), \mathcal{B}(\mathrm{supp}(\mu)), T, \mu)$ satisfies
LDT-I (resp., LDT-II, LDT-III). We obtain a surprising
result (see Theorem \ref{LDT-II'}), however, which shows that LDT-II'
is equivalent to the ergodicity.

The following example shows that ergodicity, LDT, WLDT, LDT-I, LDT-I',
LDT-II, LDT-II', LDT-III, and LDT-III' are not preserved under iterations.
This example also shows that LDT does not guarantee the weakly mixing
property. Let $\mathscr{P}_{1}$ be the set of all above listed 
properties. 

\begin{example}
Let $X=\{a_{1}, a_{2}\}$ be any two distinct points with a discrete metric
satisfying $\mu(a_{1})=\mu(a_{2})=1/2$. Define $T: X \rightarrow X$ as
$T(a_{1})=a_{2}$ and $T(a_{2})=a_{1}$. It is easy to see that, for each
$\mathrm{P}\in \mathscr{P}_{1}$, $T$ satisfies $\mathrm{P}$, but $T^{2}$ does not.
\end{example}

In \cite{WWC}, we proved that if there exists
$n\in \mathbb{N}$ such that $(T^{n}, \mu)$ satisfies LDT, then $(T, \mu)$
satisfies LDT. Similarly, it can be verified that this also holds
for all properties in $\mathscr{P}_{1}$, as summarized below.
\begin{theorem}
Let $(X, T)$ be a dynamical system, $\mu\in M(X, T)$, and $\mathrm{P}\in \mathscr{P}_{1}$.
If there exists $n\in \mathbb{N}$ such that $(T^{n}, \mu)$ satisfies $\mathrm{P}$, then $(T, \mu)$
satisfies $\mathrm{P}$.
\end{theorem}

\subsection{The central limit theorem}
Recently, Niu and Wang \cite{Niu-Wang2010} applied the central limit theorems in probability
theory to dynamical systems and found certain relations between the central limit
theorem and some chaotic properties.

An observable $\varphi$ is said to satisfy the {\it Central Limit
Theorem in the sense of Niu and Wang} for $(T, \mu)$, if there exists $\sigma>0$ such that for
every interval $A\subset \mathbb{R}$,
\[
\lim_{n\rightarrow+\infty}\mu \left(\left\{x\in X:
\frac{1}{\sqrt{n}}\sum_{i=0}^{n-1}\left(\varphi\circ T^{i}(x)- \int_{X}
\varphi \mathrm{d}\mu\right)\in A\right\}\right)=\frac{1}{\sigma
\sqrt{2\pi}}\int_{A}e^{-t^{2}/ (2\sigma^{2})}\mathrm{d}t.
\]
If every observable satisfies the central limit theorem for $(T, \mu)$ and $\mathrm{supp}(\mu)=X$,
it is said that $(T, \mu)$ satisfies the central limit theorem.

At the end of this paper, however, we will show that a dynamical system satisfying the Central Limit Theorem
in the sense of Niu and Wang actually does not exist (see Theorem \ref{CLT}).

\section{The large deviations theorem and ergodicity}\label{S-3}



In this section, some ergodic properties on dynamical systems satisfying
LDT are obtained. The new results show that LDT-I' $\Rightarrow$ LDT-II'
$\Leftrightarrow$ ergodicity.


The following lemma gives some characteristics to the ergodic measure.
\begin{lemma}\label{ergodic-char}
Let $(X, T)$ be a dynamical system and $\mu\in M(X)$. Then, the following statements
are equivalent:
\begin{enumerate}
  \item\label{er-1} $\mu$ is ergodic.
  \item\label{er-2} There exists $Y\in \mathcal{B}(X)$ with $\mu(Y)=1$ such that, for all $y\in Y$,
  $$
  \frac{1}{n}\sum_{i=0}^{n-1}\delta_{T^{i}(y)}\rightarrow \mu.
  $$
  \item\label{er-3} There exists $Y\in \mathcal{B}(X)$ with $\mu(Y)=1$ such that, for all $y\in Y$
  and all $\varphi\in C(X)$,
  $$
  \lim_{n\rightarrow\infty}\frac{1}{n}\sum_{i=0}^{n-1}\varphi \circ T^{i}(y)
  = \int_{X}\varphi \mathrm{d}\mu.
  $$
  \item\label{er-4} For any $\varphi\in C(X)$ and any $\varepsilon>0$,
  $$
  \lim_{k\rightarrow\infty}\mu\left(\bigcup_{n=k}^{+\infty}\left\{x\in X: \left|
  \frac{1}{n}\sum_{i=0}^{n-1}\varphi \circ T^{i}(x)
  -\int_{X}\varphi\mathrm{d}\mu \right|>\varepsilon \right\}\right)=0.
  $$
  \item\label{er-5} For any $\varphi\in C(X)$, $(1/ n)\sum_{i=0}^{n-1}\varphi
  \circ T^{i}(x)\rightarrow \int_{X}\varphi \mathrm{d}\mu$ $\mathrm{a.e.}$
\end{enumerate}
\end{lemma}
\begin{proof}
(\ref{er-2}) $\Leftrightarrow$ (\ref{er-3}) follows from Lemma \ref{lem:2-1}. 
In view of \cite[Theorem 6.14]{Walters1982}, it suffices to check that (\ref{er-3}) $\Leftrightarrow$
(\ref{er-4}), and that $\mu \in M(X, T)$ under the assumption of (\ref{er-3}).

(\ref{er-3}) $\Rightarrow$ (\ref{er-4}).
Given any fixed $\varphi\in C(X)$, one has
\begin{equation}\label{e-4-1}
\begin{split}
X\setminus Y\supset \mathrm{Div}(\varphi)&\ := \
\left\{x\in X: \frac{1}{n}\sum_{i=0}^{n-1}\varphi \circ T^{i}(x)\nrightarrow \int_{X}
\varphi \mathrm{d}\mu\right\}\\
&\ =\ \bigcup_{n=1}^{+\infty}\bigcap_{m=1}^{+\infty}\bigcup_{k=m}^{+\infty}
\left\{x\in X: \left|\frac{1}{k}\sum_{i=0}^{k-1}\varphi \circ T^{i}(x)- \int_{X}
\varphi \mathrm{d}\mu\right|> \frac{1}{n}\right\}.
\end{split}
\end{equation}
Then, $\mu(\mathrm{Div}(\varphi))=0$. So, for any
$n\in \mathbb{N}$,
$$
\lim_{m\rightarrow+\infty}\left(\bigcup_{k=m}^{+\infty}
\left\{x\in X: \left|\frac{1}{k}\sum_{i=0}^{k-1}\varphi \circ T^{i}(x)- \int_{X}
\varphi \mathrm{d}\mu\right|> \frac{1}{n}\right\}\right)=0.
$$
This implies that, for any $\varepsilon>0$,
$$
\lim_{m\rightarrow+\infty}\left(\bigcup_{k=m}^{+\infty}
\left\{x\in X: \left|\frac{1}{k}\sum_{i=0}^{k-1}\varphi \circ T^{i}(x)- \int_{X}
\varphi \mathrm{d}\mu\right|> \varepsilon\right\}\right)=0.
$$

(\ref{er-4}) $\Rightarrow$ (\ref{er-3}). Given any fixed $\varphi\in C(X)$,
condition (\ref{er-4}) implies that, for any $n\in \mathbb{N}$,
\begin{eqnarray*}
&&\quad \mu\left(\bigcap_{m=1}^{+\infty}\bigcup_{k=m}^{+\infty}
\left\{x\in X: \left|\frac{1}{k}\sum_{i=0}^{k-1}\varphi \circ T^{i}(x)- \int_{X}
\varphi \mathrm{d}\mu\right|> \frac{1}{n}\right\}\right)\\
&&=
\lim_{m\rightarrow\infty}\mu\left(\bigcup_{k=m}^{+\infty}
\left\{x\in X: \left|\frac{1}{k}\sum_{i=0}^{k-1}\varphi \circ T^{i}(x)- \int_{X}
\varphi \mathrm{d}\mu\right|> \frac{1}{n}\right\}\right)=0.
\end{eqnarray*}
Combining this with
\begin{eqnarray*}
\mathrm{Con}(\varphi)& := &
\left\{x\in X: \frac{1}{n}\sum_{i=0}^{n-1}\varphi \circ T^{i}(x)\rightarrow \int_{X}
\varphi \mathrm{d}\mu\right\}\\
& = & \bigcap_{n=1}^{+\infty}\bigcup_{m=1}^{+\infty}\bigcap_{k=m}^{+\infty}
\left\{x\in X: \left|\frac{1}{k}\sum_{i=0}^{k-1}\varphi \circ T^{i}(x)- \int_{X}
\varphi \mathrm{d}\mu\right|< \frac{1}{n}\right\},
\end{eqnarray*}
it follows that $\mu(\mathrm{Con}(\varphi))=1$.
Choose a countable dense subset $\{\varphi_{k}\}_{k=1}^{+\infty}$ of $C(X)$
and take $Y=\bigcap_{k=1}^{+\infty}\mathrm{Con}(\varphi_{k})$. Then, $\mu(Y)=1$
and, for any $y\in Y$ and any $k\in \mathbb{N}$,
$$
\lim_{n\rightarrow +\infty}\frac{1}{n}\sum_{i=0}^{n-1}\varphi_{k} \circ T^{i}(y)= \int_{X}\varphi_{k}
  \mathrm{d}\mu.
$$
The result follows from approximating a given $\varphi\in C(X)$ by members of $\{\varphi_{k}\}_{k=1}^{+\infty}$.

Applying \eqref{e-4-1}, it is easy to see (\ref{er-4}) $\Leftrightarrow$ (\ref{er-5}).

Finally, according to the proof of the Krylov-Bogolioubov Theorem (see \cite[Theorem 6.9]{Walters1982}),
it is easy to see that under the assumption of
(\ref{er-3}), $\mu$ is $T$-invariant, i.e., $\mu\in M(X, T)$. Its
proof is included here for completeness.

In fact, for any given $\varphi\in C(X)$ and $y\in Y$, noting that $\varphi \circ T\in C(X)$,
condition (\ref{er-3}) implies that $(1/ n)\sum_{i=0}^{n-1}\varphi \circ T^{i}(y)
\rightarrow \int_{X}\varphi \mathrm{d}\mu$ and $(1/ n)\sum_{i=0}^{n-1}(\varphi\circ T) \circ T^{i}(y)
\rightarrow \int_{X}\varphi\circ T \mathrm{d}\mu$. Combining this with
$|(1/ n)\sum_{i=0}^{n-1}\varphi \circ T^{i}(y)-
(1/n)\sum_{i=0}^{n-1}(\varphi\circ T) \circ T^{i}(y)|\leq 2\|\varphi\|/ n$,
it follows that
$$
\int_{X}\varphi\mathrm{d}\mu-
\int_{X}\varphi\circ T \mathrm{d}\mu=
\lim_{n\rightarrow +\infty}\frac{1}{n}\sum_{i=0}^{n-1}\varphi \circ T^{i}(y)-
\lim_{n\rightarrow +\infty}\frac{1}{n}\sum_{i=0}^{n-1}(\varphi\circ T) \circ T^{i}(y)=0.
$$
The result is now implied  by \cite[Theorem 6.8]{Walters1982}.
\end{proof}

\begin{remark}\label{remark-*}
\begin{enumerate}
  \item It is noticeable that \cite[Theorem 6.14]{Walters1982} shows that
  Lemma \ref{ergodic-char} (\ref{er-1}) is equivalent to
  Lemma \ref{ergodic-char} (\ref{er-2}), under the hypothesis of
  $\mu\in M(X, T)$.
  \item By Lemma \ref{ergodic-char} (\ref{er-3}), it is easy to see
  that $\mu\in M(X, T)$ is ergodic if and only if $\mu\left(\left\{x\in X:
  \lim_{n\rightarrow+\infty}(1/ n)\sum_{i=0}^{n-1}\varphi \circ T^{i}(x)=\int_{X}\varphi\mathrm{d}\mu,
  \ \forall \varphi\in C(X)\right\}\right)=1$.
\end{enumerate}
\end{remark}

\begin{theorem}\label{Th 3-1}
Let $(X, T)$ be a dynamical system and $\mu\in M(X)$. If $(T, \mu)$ satisfies
LDT-I', then $\mu$ is ergodic.
\end{theorem}
\begin{proof}
The result follows from
Lemma \ref{ergodic-char} and the fact that
\begin{eqnarray*}
&&\quad \mu\left(\bigcup_{n=k}^{+\infty}\left\{x\in X: \left|
  \frac{1}{n}\sum_{i=0}^{n-1}\varphi \circ T^{i}(x)
  -\int_{X}\varphi\mathrm{d}\mu \right|>\varepsilon \right\}\right)\\
&&  \leq \sum_{n=k}^{+\infty}\mu\left(\left\{x\in X: \left|
  \frac{1}{n}\sum_{i=0}^{n-1}\varphi \circ T^{i}(x)
  -\int_{X}\varphi\mathrm{d}\mu \right|>\varepsilon\right\}\right).
\end{eqnarray*}
\end{proof}

Because LDT is stronger than
LDT-I', Theorem \ref{Th 3-1} immediately generates the
following corollaries.
\begin{corollary}\label{Corollary **}
If $(T, \mu)$ satisfies WLDT, then $\mu$
ergodic.
\end{corollary}
\begin{corollary}\label{Corollary **}
If $(T, \mu)$ satisfies LDT or LDT-I, then $\mu$
is an ergodic measure with a full support.
\end{corollary}

\begin{corollary}\cite[Theorem 3.1]{Gu2007}\label{Gu}
If $(T, \mu)$ satisfies LDT, then
$T$ is topologically ergodic.
\end{corollary}
\begin{proof}
Since every ergodic dynamical system with full support is $\mathscr{F}_{s}$-transitive,
this follows by Corollary \ref{Corollary **}.
\end{proof}

\begin{theorem}\label{ue-ldt}
Let $(X, T)$ be a uniquely ergodic dynamical system and $\mu\in E(X, T)$. Then,
\begin{enumerate}
\item\label{ue-1} $(T, \mu)$ satisfies WLDT.
\item $(T|_{\mathrm{supp}(\mu)}, \mu)$ satisfies LDT.
\end{enumerate}
\end{theorem}
\begin{proof}
The unique ergodicity of $(X, T)$, together with
Lemma \ref{Unique ergodic}, implies that for every $\varphi\in C(X)$,
$(1/ n)\sum_{i=0}^{n-1}\varphi\circ T^{i}(x)$ converges uniformly
  to $\int_{X}\varphi\mathrm{d}\mu$. This means that for any
  $\varepsilon>0$, there exists $N\in \mathbb{N}$ such that for
  all $x\in X$, $\left|(1/ n)\sum_{i=0}^{n-1}\varphi\circ T^{i}(x)
  -\int_{X}\varphi\mathrm{d}\mu\right|<\varepsilon$ holds for all
  $n\geq N$. So, $(T, \mu)$ satisfies WLDT.

(2) This follows by (\ref{ue-1}) and $\mathrm{supp}(\mu|_{\mathcal{B}(\mathrm{supp}(\mu))})=\mathrm{supp}(\mu)$.
\end{proof}

Theorem \ref{ue-ldt} with Lemma \ref{ergodic-support} leads
 to the following result.

\begin{corollary}\label{Corollary 3-3-3}
Let $(X, T)$ be a uniquely ergodic dynamical system and $\mu\in M(X)$. Then,
the following statements are equivalent:
\begin{enumerate}
  \item $(T, \mu)$ satisfies LDT.
  \item $(T, \mu)$ satisfies LDT-I.
  \item $(T, \mu)$ satisfies LDT-II.
  \item $T$ is minimal and $\mu\in M(X, T)$.
\end{enumerate}
\end{corollary}

\begin{remark}\label{Remark-2}
  Corollary \ref{Corollary 3-3-3} indicates that a nontrivial dynamical
  system satisfying LDT, which is sensitive or equicontinuous,
  indeed exists.
\end{remark}

The following Proposition \ref{LDT-II} shows that the probability measure $\mu$
of the pair $(T, \mu)$ satisfying LDT-II' is $T$-invariant.
Theorem \ref{LDT-II'} gives an improved characteristic of LDT-II'
which indicates that LDT-II' is equivalent to ergodicity.

\begin{proposition}\label{LDT-II}
Let $(X, T)$ be a dynamical system and $\mu\in M(X)$. If $(T, \mu)$ satisfies
LDT-II', then $\mu$ is $T$-invariant.
\end{proposition}
\begin{proof}
By \cite[Theorem 6.8]{Walters1982}, it suffices to check that, for any $\varphi\in C(X)$,
$\int_{X} \varphi \mathrm{d}\mu=\int_{X} \varphi\circ T \mathrm{d}\mu$.
Since $\varphi \circ T\in C(X)$ and $(T, \mu)$ satisfies LDT-II',
it follows that, for any $M\in \mathbb{N}$,
$$
\lim_{n\rightarrow +\infty}\mu\left(\left\{x\in X: \left|
\frac{1}{n}\sum_{i=0}^{n-1}\varphi \circ T^{i}(x)
-\int_{X}\varphi\mathrm{d}\mu \right|\leq \frac{1}{2M}\right\}\right)=1
$$
and
$$
\lim_{n\rightarrow +\infty}\mu\left(\left\{x\in X: \left|
\frac{1}{n}\sum_{i=0}^{n-1}(\varphi\circ T) \circ T^{i}(x)
-\int_{X}\varphi \circ T \mathrm{d}\mu \right|\leq \frac{1}{2M}\right\}\right)=1.
$$
Thus, there exists $N\in \mathbb{N}$ such that, for any
$n\geq N$, there exists $x_{n}\in X$ satisfying
$|(1/ n)\sum_{i=0}^{n-1}\varphi \circ T^{i}(x_{n})
-\int_{X}\varphi\mathrm{d}\mu |\leq 1/ 2M$ and $|
(1/ n)\sum_{i=0}^{n-1}(\varphi\circ T) \circ T^{i}(x_{n})
-\int_{X}\varphi \circ T \mathrm{d}\mu |
\leq 1/ 2M$. This implies that
\begin{eqnarray*}
&&\left|\int_{X} \varphi \mathrm{d}\mu-\int_{X} \varphi\circ T \mathrm{d}\mu\right|\\
& \leq & \left|\int_{X} \varphi \mathrm{d}\mu-\frac{1}{n}\sum_{i=0}^{n-1}\varphi \circ T^{i}(x_{n})\right|
+\left|\frac{1}{n}\sum_{i=0}^{n-1}\varphi \circ T^{i}(x_{n})-\frac{1}{n}\sum_{i=0}^{n-1}(\varphi\circ T) \circ T^{i}(x_{n})\right|\\
&&  + \left|\frac{1}{n}\sum_{i=0}^{n-1}(\varphi\circ T) \circ T^{i}(x_{n})-\int_{X} \varphi\circ T \mathrm{d}\mu\right|\\
& \leq & \frac{1}{2M} +\frac{2\|\varphi\|}{n}+ \frac{1}{2M}=\frac{1}{M} +\frac{2\|\varphi\|}{n}.
\end{eqnarray*}
So, $\int_{X} \varphi \mathrm{d}\mu=\int_{X} \varphi\circ T \mathrm{d}\mu$.
\end{proof}

\begin{theorem}\label{LDT-II'}
Let $(X, T)$ be a dynamical system and $\mu\in M(X)$. Then, $(T, \mu)$
satisfies LDT-II' if and only if $\mu$ is ergodic.
\end{theorem}
\begin{proof}
The sufficiency follows immediately from Lemma \ref{ergodic-char}.

To prove the
necessity, by Lemma \ref{ergodic-char}, it suffices to check that, for any $\varphi\in C(X)$,
$(1/ n)\sum_{i=0}^{n-1}\varphi
  \circ T^{i}(x)\rightarrow \int_{X}\varphi \mathrm{d}\mu$ a.e.

\begin{claim}\label{Claim 1}
For any $\varphi\in C(X)$, there exists an increasing sequence $\{L_{k}\}_{k=1}^{+\infty}\subset \mathbb{N}$
such that $(1/ L_{k})\sum_{i=0}^{L_{k}-1}\varphi
  \circ T^{i}(x)\rightarrow \int_{X}\varphi \mathrm{d}\mu$ $\mathrm{a.e.}$
\end{claim}

Given any $N\in \mathbb{N}$, since $(T, \mu)$ satisfies LDT-II', one has
$$
\lim_{n \rightarrow +\infty}\mu\left(\Sigma(n, 1/N):=\left\{x\in X: \left|
\frac{1}{n}\sum_{i=0}^{n-1}\varphi \circ T^{i}(x)
-\int_{X}\varphi\mathrm{d}\mu \right|>\frac{1}{N}\right\}\right)=0.
$$
In particular, for $N=1$, there exists an increasing sequence
$\{L^{(1)}_{k}\}_{k=1}^{+\infty}$ such that, for any $k\in \mathbb{N}$,
$\mu(\Sigma(L^{(1)}_{k}, 1))\leq 1/ 2^{k}$. Consider
the sequence of numbers $\{\mu (\Sigma(L^{(1)}_{k}, 1/ 2) ) \}_{k=1}^{+\infty}$,
which converges to zero by the fact that $(T, \mu)$ satisfies
LDT-II' and so has a subsequence
$\{\mu (\Sigma(L^{(2)}_{k}, 1/ 2))\}_{k=1}^{+\infty}$ such that
for any $k\in \mathbb{N}$, $\mu (\Sigma(L^{(2)}_{k}, 1/2) )\leq 1/ 2^{k}$.
Clearly, $\sum_{k=1}^{+\infty}\mu (\Sigma(L^{(2)}_{k}, 1) )\leq \sum_{k=1}^{+\infty}
1/ 2^{k}=1$. Proceed in this process and, for each $N\in \mathbb{N}$, obtain a subsequence
$\{L^{(N)}_{k}\}_{k=1}^{+\infty}\subset \mathbb{N}$ such that
$\{L^{(N)}_{k}\}_{k=1}^{+\infty}\subset \{L^{(N-1)}_{k}\}_{k=1}^{+\infty}
\subset \cdots \subset \{L^{(2)}_{k}\}_{k=1}^{+\infty}\subset \{L^{(1)}_{k}\}_{k=1}^{+\infty}$,
so that $\sum_{k=1}^{+\infty}\mu (\Sigma(L^{(j)}_{k}, 1/ j) )\leq 1$ for
$j=1, \ldots, N$. Choose the diagonal $\{L_{k}=L^{(k)}_{k}\}_{k=1}^{+\infty}$. It is easy
to see that $\sum_{k=1}^{+\infty}\mu\left(\Sigma(L_{k}, 1/ j)\right)\leq 1$ holds for all
$j=1, 2, \ldots$. Combining this with
\begin{eqnarray*}
Q & := &
\left\{x\in X: \frac{1}{L_{k}}\sum_{i=0}^{L_{k}-1}\varphi \circ T^{i}(x)\nrightarrow \int_{X}
\varphi \mathrm{d}\mu\right\}\\
& = &\bigcup_{n=1}^{+\infty}\bigcap_{m=1}^{+\infty}\bigcup_{k=m}^{+\infty}
\left\{x\in X: \left|\frac{1}{L_{k}}\sum_{i=0}^{L_{k}-1}\varphi \circ T^{i}(x)- \int_{X}
\varphi \mathrm{d}\mu\right|> \frac{1}{n}\right\},
\end{eqnarray*}
it follows that $\mu(Q)=0$. The proof of Claim \ref{Claim 1} is thus completed.

\begin{claim}
For any $\varphi\in C(X)$, $(1/ n)\sum_{i=0}^{n-1}\varphi
  \circ T^{i}(x)\rightarrow \int_{X}\varphi \mathrm{d}\mu$ $\mathrm{a.e.}$
\end{claim}

By the Birkhoff Ergodic Theorem, Proposition \ref{LDT-II}, and Claim \ref{Claim 1},
there exist $\varphi^{*}\in L^{1}(\mu)$ (as $\varphi\in C(X)\subset
L^{1}(\mu)$) and an increasing
sequence $\{L_{k}\}_{k=1}^{+\infty}$
such that $(1/ n)\sum_{i=0}^{n-1}\varphi
  \circ T^{i}(x)\rightarrow \varphi^{*}(x)$ a.e. and
  $(1/ L_{k})\sum_{i=0}^{L_{k}-1}\varphi
  \circ T^{i}(x)\rightarrow \int_{X}\varphi \mathrm{d}\mu$ a.e.
  Take $\Sigma=\{x\in X: (1/ n)\sum_{i=0}^{n-1}\varphi
  \circ T^{i}(x)\rightarrow \varphi^{*}(x)\} \bigcap
  \{x\in X: (1/ L_{k})\sum_{i=0}^{L_{k}-1}\varphi
  \circ T^{i}(x)\rightarrow \int_{X}\varphi \mathrm{d}\mu\}$.
  Clearly, $\mu(\Sigma)=1$ and, for any $x\in \Sigma$,
$$
\varphi^{*}(x)=\lim_{n\rightarrow+\infty}\frac{1}{n}\sum_{i=0}^{n-1}\varphi
  \circ T^{i}(x)=\lim_{k\rightarrow+\infty}
  \frac{1}{L_{k}}\sum_{i=0}^{L_{k}-1}\varphi
  \circ T^{i}(x)=\int_{X}\varphi \mathrm{d}\mu.
$$
This implies that $(1/ n)\sum_{i=0}^{n-1}\varphi
  \circ T^{i}(x)\rightarrow \int_{X}\varphi \mathrm{d}\mu$ a.e.
\end{proof}

In the proof of Theorem \ref{LDT-II'}, the diagonalization procedure works.
We will give another proof of this result in Section \ref{S-6}
(see Lemma \ref{Lemma 6.1}).



\begin{corollary}\label{Corollary 3.2}
Let $(X, T)$ be a dynamical system and $\mu\in M(X)$. If
$(T, \mu)$ satisfies LDT, LDT-I or LDT-II,
then $(X, T)$ is an $E$-system. In particular, $T$ is $\mathscr{F}_{s}$-transitive.
\end{corollary}
\begin{proof}
The result yields by Theorem \ref{LDT-II'}.
\end{proof}



\section{The large deviations theorem and sensitivity}

An interesting question about a dynamical system is when its orbits from nearby
points start to separate after finite steps. This is also one of the most important features
depicting the chaoticity of a system. This notion, referred to as
the ``butterfly effect'', has been widely studied and is termed the {\it sensitive
dependence on initial conditions} (briefly, sensitivity),
introduced by Auslander and Yorke~\cite{AY80} and popularized
by Devaney~\cite{Devaney1989}.

In \cite{Huang-Ye2002b}, Huang and Ye proved that an almost equicontinuous
$\mathscr{F}_{s}$-transitive system is minimal and equicontinuous. This, together with
the Auslander-Yorke dichotomy Theorem (also see \cite[Theorem 3.1]{Akin-Kolyada2003})
and \cite[Corollary 1]{Moothathu2007}, which states that for an $\mathscr{F}_{s}$-transitive
system, sensitivity implies syndetical sensitivity, implies that every
$\mathscr{F}_{s}$-transitive system is either syndetically sensitive or both minimal and equicontinuous.
These with Corollary \ref{Corollary 3.2} lead to the following results.

\begin{theorem}\label{weak-mix}
If $(X, T)$ is strongly topologically ergodic, then $T$ is weakly mixing. In particular,
$T$ is sensitive.
\end{theorem}
\begin{proof}
It follows from the fact that the weakly mixing property is equivalent to the $\mathscr{F}_{t}$-transitivity
and $\mathscr{F}_{t}\supset \mathscr{F}_{ud1}$.
\end{proof}

\begin{theorem}\label{sen-Th-1}
If $(T, \mu)$ satisfies LDT-II and
$T$ is sensitive, then $T$ is
syndetically sensitive.
\end{theorem}

Theorem \ref{sen-Th-1} implies \cite[Theorem 4.1]{Gu2007} immediately.

\begin{corollary}\label{sen-Th}
If $(T, \mu)$ satisfies LDT-II and
$T$ is strongly topologically ergodic, then $T$ is
syndetically sensitive.
\end{corollary}

\begin{corollary}\cite[Theorem 3.1]{Li2013}\label{Li}
If $(T, \mu)$ satisfies LDT and
$T$ is strongly topologically ergodic, then $T$ is
ergodically sensitive.
\end{corollary}

\begin{corollary}\label{Corollary 3.3}
If $(T, \mu)$ satisfies LDT-II,
then $T$ is either syndetically sensitive or both minimal and equicontinuous.
\end{corollary}

For the $\mathscr{F}$-transitivity, we have the following result.
Note that the proof is similar to \cite[Theorem 4.6]{Huang-Ye2002b},
but for completeness, a proof is provided here.

\begin{theorem}\label{Thm-De}
Let $(X, T)$ be a dynamical system and let $\mathscr{F}$ be a Furstenberg family.
If $(X, T)$ is $\mathscr{F}$-transitive and almost equicontinuous, then
for any $x\in \mathrm{Trans}(T)$ and any neighbourhood $U$ of $x$,
$N(x, U)\in \Delta(\mathscr{F})$.
\end{theorem}
\begin{proof}
Fix any $x\in \mathrm{Trans}(T)$ and any $\varepsilon>0$. Noting that $x$ is an equicontinuous
point, it follows that there exists $0<\varepsilon_{1}<\varepsilon/ 2$ such that, for any $y\in B(x, \varepsilon_{1})$
and any $n\in \mathbb{Z}^{+}$, $d(T^{n}(x), T^{n}(y))<\varepsilon/ 2$. For any $n\in N(B(x, \varepsilon_{1}), B(x, \varepsilon_{1}))$,
there exists $y\in B(x, \varepsilon_{1})$ such that $T^{n}(y)\in B(x, \varepsilon_{1})$. Thus, $d(x, T^{n}(x))
\leq d(x, y)+ d(T^{n}(x), T^{n}(y))<\varepsilon$, i.e., $N(B(x, \varepsilon_{1}), B(x, \varepsilon_{1}))
\subset N(x, B(x, \varepsilon))\in \mathscr{F}$. This, together with \cite[Proposition 2.2]{YZ},
implies that for any neighbourhood $U$ of $x$, $N(U, U)=N(x, U)-N(x, U)\in \Delta(\mathscr{F})$.
Repeating the proof above, one obtains that $N(x, U)\in \Delta(\mathscr{F})$.
\end{proof}

Now, consider the non-sensitive dynamical systems satisfying LDT-II'.


\begin{theorem}\label{Theorem 3-4}
Let $(X, T)$ be a dynamical system and $(T, \mu)$ satisfy
LDT-II'. If $(X, T)$ is not sensitive, then
\begin{enumerate}
\item\label{3-5-1} $(\mathrm{supp}(\mu), T)$ is minimal and equicontinuous.
\item\label{3-5-2} For every observable $\varphi\in C(\mathrm{supp}(\mu))$,
$$
\lim_{n\rightarrow +\infty}\frac{1}{n}\sum_{i=0}^{n-1}\varphi \circ T^{i}(x)=\int_{\mathrm{supp}(\mu)}
\varphi\mathrm{d}\mu, \quad \forall x\in \mathrm{supp}(\mu).
$$
\item\label{3-5-3} $(\mathrm{supp}(\mu), T)$ is uniquely ergodic.
\end{enumerate}
\end{theorem}

\begin{proof}
(1) It follows directly from Corollary \ref{Corollary 3.2} and Theorem \ref{Thm-De}.

(2) Suppose, on the contrary, that there exist an observable $\varphi\in C(\mathrm{supp}(\mu))$
and $x\in \mathrm{supp}(\mu)$ such that $(1/ n)\sum_{i=0}^{n-1}\varphi \circ T^{i}(x)\nrightarrow \int_{\mathrm{supp}(\mu)}
\varphi\mathrm{d}\mu$, i.e.,
$$
\xi:=\limsup_{n\rightarrow +\infty}\frac{1}{n}\sum_{i=0}^{n-1}\varphi \circ T^{i}(x)\neq \int_{\mathrm{supp}(\mu)}
\varphi\mathrm{d}\mu,
$$
or
$$
\eta:=\liminf_{n\rightarrow +\infty}\frac{1}{n}\sum_{i=0}^{n-1}\varphi\circ T^{i}(x)
\neq \int_{\mathrm{supp}(\mu)} \varphi\mathrm{d}\mu.
$$
Without loss of generality, assume that $\xi>\int_{\mathrm{supp}(\mu)} \varphi\mathrm{d}\mu$,
because the rest cases can be verified similarly. Then, there exists an
increasing sequence $\{n_{k}\}_{k=1}^{+\infty}\subset\mathbb{N}$ such that, for any
$k\in \mathbb{N}$,
\begin{equation}\label{e-1}
\frac{1}{n_{k}}\sum_{i=0}^{n_{k}-1}\varphi\circ T^{i}(x)> \frac{1}{4}\left(3\xi +\int_{\mathrm{supp}(\mu)} \varphi\mathrm{d}\mu\right).
\end{equation}
Since $\varphi$ is uniformly continuous, there exists $0<\delta_{1}<(\xi-\int_{\mathrm{supp}(\mu)} \varphi\mathrm{d}\mu) / 4$
such that, for any $x_{1}, x_{2}\in X$ with $d(x_{1}, x_{2})<\delta_{1}$,
\begin{equation}\label{e-2}
|\varphi(x_{1})-\varphi(x_{2})|<\frac{1}{4}\left(\xi-\int_{\mathrm{supp}(\mu)} \varphi\mathrm{d}\mu\right).
\end{equation}
The equicontinuity of $T$ implies that there exists $0<\delta<\delta_{1}$ such that,
for any $y\in B(x, \delta):=\left\{y\in \mathrm{supp}(\mu): d(x, y)<\delta\right\}$ and any
$n\in \mathbb{Z}^{+}$, $d(T^{n}(x), T^{n}(y))<\delta_{1}$.
Combining this with (\ref{e-1}) and (\ref{e-2}),
for any $y\in B(x, \delta)$ and any $k\in \mathbb{N}$, one has
\begin{eqnarray*}
\frac{1}{n_{k}}\sum_{i=0}^{n_{k}-1}\varphi\circ T^{i}(y) & = &
\frac{1}{n_{k}}\sum_{i=0}^{n_{k}-1}\varphi\circ T^{i}(x)-\frac{1}{n_{k}}\sum_{i=0}^{n_{k}-1}
\left[\varphi \circ T^{i}(x)-\varphi\circ T^{i}(y)\right]\\
& \geq & \frac{1}{n_{k}}\sum_{i=0}^{n_{k}-1}\varphi\circ T^{i}(x)-
\frac{1}{n_{k}}\sum_{i=0}^{n_{k}-1}\left|\varphi\circ T^{i}(x)-\varphi\circ T^{i}(y)\right|\\
& > & \frac{1}{4}\left(3\xi +\int_{\mathrm{supp}(\mu)} \varphi\mathrm{d}\mu\right)-
\frac{1}{4}\left(\xi-\int_{\mathrm{supp}(\mu)} \varphi\mathrm{d}\mu\right)\\
& = & \frac{1}{2}\left(\xi+\int_{\mathrm{supp}(\mu)} \varphi\mathrm{d}\mu\right),
\end{eqnarray*}
i.e.,
$$
\frac{1}{n_{k}}\sum_{i=0}^{n_{k}-1}\varphi\circ T^{i}(y)-\int_{\mathrm{supp}(\mu)} \varphi\mathrm{d}\mu
>\frac{1}{2}\left(\xi-\int_{\mathrm{supp}(\mu)} \varphi\mathrm{d}\mu \right)>0.
$$
This implies that, for any $k\in \mathbb{N}$,
$$
\left\{z\in \mathrm{supp}(\mu): \left|\frac{1}{n_{k}}\sum_{i=0}^{n_{k}-1}\varphi\circ T^{i}(z)-\int_{\mathrm{supp}(\mu)} \varphi\mathrm{d}\mu\right|
>\frac{1}{2}\left(\xi-\int_{\mathrm{supp}(\mu)} \varphi\mathrm{d}\mu\right)\right\}\supset B(x, \delta),
$$
which is a contradiction since $(T|_{\mathrm{supp}(\mu)}, \mu)$ satisfies
LDT-II and $\mu(B(x, \delta))>0$.

(3) Based on Lemma \ref{Unique ergodic}, Proposition \ref{LDT-II},
and (\ref{3-5-2}), this holds trivially.
\end{proof}

\begin{corollary}
Let $T$ be a measure-preserving transformation of a probability space $(X, \mathcal{B}(X, \mu)$.
Then, $(X, T)$ is a non-sensitive ergodic system with full support
if and only if $T$ is minimal, equicontinuous,
and uniquely ergodic.
\end{corollary}
\begin{proof}
It follows immediately from Corollary \ref{Corollary 3-3}, Theorem \ref{LDT-II'}, and
Theorem \ref{Theorem 3-4}.
\end{proof}

\section{The large deviations theorem of type III and transitivity}\label{S-6}

This section studies the transitivity of dynamical
systems satisfying LDT-III.

\begin{theorem}\label{Lemma 1}
Let $(X, T)$ be a dynamical system and $\mu\in M(X)$.
If $(T, \mu)$ satisfies LDT-III', then:
\begin{enumerate}
  \item\label{LDT-III-1} For any nonempty open subset $U\subset X$ with $\mu(U)>0$, there exists
  $x\in U$ such that $\overline{d}(N(x, U))>0$.
  \item For any nonempty open subsets $U, V\subset X$ with $\mu(U)
  \mu(V)>0$, $N(U, V)\neq \O$.
\end{enumerate}
\end{theorem}
\begin{proof} (1) Suppose, on
the contrary, that (\ref{LDT-III-1}) does not hold. Then, there exists a
nonempty open set $U\subset X$ with $\mu(U)>0$ such that, for any $x\in U$, $\overline{d}(N(x, U))=0$,
i.e.,
$$
\lim_{n\rightarrow +\infty}\frac{1}{n}\left|\left\{0\leq i< n:
T^{i}(x)\in U\right\}\right|=0.
$$

Take a nonempty open subset $V \subset U$ satisfying
$\overline{V}\subset U$. Clearly, both $\overline{V}$ and
$X\setminus U$ are nonempty closed subsets of $X$,
and $\overline{V} \cap (X\setminus U)=\O$. Applying Urysohn's Lemma,
there exists a continuous function $\varphi: X\rightarrow
[0, 1]$ such that
\[
\varphi(x)=\left\{\begin{array}{cc}
1, & {\rm } x\in \overline{V}, \\\
0, & {\rm } x\in X\setminus U.
\end{array}
\right.
\]
Clearly, $\varphi$ is an observable.
Choose $\varepsilon=\frac{1}{3}\int_{X} \varphi \mathrm{d}\mu>0$ and
set
$$
D_{n}^{(\varepsilon)}=\left\{x\in X: \left|
\frac{1}{n}\sum_{i=0}^{n-1}\varphi\circ T^{i}(x)-\int_{X} \varphi \mathrm{d}\mu
\right|>\varepsilon\right\}, \ n=1, 2, \ldots.
$$
Since the pair $(T, \mu)$ satisfies LDT-III', one has
$$
\liminf_{n\rightarrow +\infty}\mu\left(D_{n}^{(\varepsilon)}\right)=0.
$$
This implies that there exists an increasing sequence $\left\{N_{k}\right\}_{k=1}^{+\infty}
\subset \mathbb{N}$
such that, for any $k\in \mathbb{N}$,
\begin{equation}\label{1}
\mu\left(D_{N_{k}}^{(\varepsilon)}\right)\leq \frac{\mu ({U})}{2^{k+1}}.
\end{equation}
It follows from the choice of $\varphi$ that, for any $x\in U$,
$$
\left|\frac{1}{n}\sum_{i=0}^{n-1}\varphi \circ T^{i}(x)\right|\leq
\frac{1}{n}\left|\left\{0\leq i< n:
T^{i}(x)\in U\right\}\right|\rightarrow 0, \
(n\rightarrow +\infty).
$$
Thus,
$$
\lim_{n\rightarrow +\infty}
\left|\frac{1}{n}\sum_{i=0}^{n-1}\varphi\circ T^{i}(x)-\int_{X} \varphi
\mathrm{d}\mu\right|=3\varepsilon.
$$
So, for any $x\in U$, there exists $M_{x}\in \mathbb{Z}^{+}$ such that, for any
$n\geq M_{x}$, $x\in D_{n}^{(\varepsilon)}$. This implies that, for any increasing sequence
$\left\{n_{k}\right\}_{k=1}^{+\infty}\subset \mathbb{N}$, $U \subset \bigcup_{k=1}^{+\infty}D_{n_{k}}^{(\varepsilon)}$.
Combining this with \eqref{1} yields that
$$
\mu(U)\leq \mu\left(\bigcup_{k=1}^{+\infty}D_{N_{k}}^{(\varepsilon)}\right)
\leq \sum_{k=1}^{+\infty}\mu \left(D_{N_{k}}^{(\varepsilon)}\right)
\leq \sum_{k=1}^{+\infty}\frac{\mu(U)}{2^{k+1}}=\frac{1}{2}\mu (U),
$$
which is a contradiction as $\mu (U)>0$. 

(2) Suppose that there exist nonempty open subsets $U, V\subset X$  with $\mu(U)
\mu(V)>0$ such that $N(U, V)=\O$, i.e., for any $n\in \mathbb{Z}^{+}$, $T^{n}(U)\cap V=\O$.
Take $F_{1}=\overline{\cup_{n=0}^{+\infty}T^{n}(U)}$. Clearly, $F_{1}\subset X\setminus V$.
Since the probability measure $\mu$ is regular and $\mu(V)>0$, there exists a closed subset
$F_{2}\subset V$ such that $\mu(F_{2})>0$. Noting that $F_{1}\cap F_{2}=\O$, by Urysohn's Lemma,
there exists a continuous function $\varphi: X\rightarrow
[0, 1]$ such that
\[
\varphi(x)=\left\{\begin{array}{cc}
1, & {\rm } x\in F_{2}, \\\
0, & {\rm } x\in F_{1}.
\end{array}
\right.
\]
It is easy to see that for any $x\in U$ and any $i\in \mathbb{Z}^{+}$,
$\varphi\circ T^{i}(x)=0$. So,
$$
U\subset \bigcap_{n=1}^{+\infty}\left\{x\in X: \left|\frac{1}{n}\sum_{i=0}^{n-1}\varphi \circ T^{i}(x)
-\int_{X}\varphi \mathrm{d}\mu\right|>\frac{\mu(F_{2})}{2}\right\},
$$
which is a contradiction as $(T, \mu)$ satisfies LDT-III' and $\mu (U)>0$.
\end{proof}

Although we do not know if the probability measure $\mu$ is $T$-invariant when
$(T, \mu)$ satisfies LDT-III or LDT-III', Corollary
\ref{Corollary ***} below indicates that such a system admits a $T$-invariant probability
measure $\nu\in M(X, T)$ such that $\mathrm{supp}(\nu)=\mathrm{supp}(\mu)$.
\begin{theorem}\label{Theorem 1}
Let $(X, T)$ be a dynamical system and $\mu\in M(X)$.
If $(T, \mu)$ satisfies LDT-III, then
$(X, T)$ is an $E$-system. In particular, $T$ is $\mathscr{F}_{s}$-transitive.
\end{theorem}
\begin{proof}
Based on Theorem~\ref{Lemma 1} and \cite[Corollary 2.2]{Yin-Zhou2012}, this
holds trivially.
\end{proof}

\begin{corollary}\label{Corollary ***}
Let $(X, T)$ be a dynamical system and $\mu\in M(X)$.
If $(T, \mu)$ satisfies LDT-III', then
$(\mathrm{supp}(\mu), T)$ is an $E$-system.
\end{corollary}

\begin{lemma}\label{Lemma 6.1}
Let $(X, T)$ be a dynamical system and $\mu\in M(X, T)$.
If $(T, \mu)$ satisfies LDT-III', then $\mu$ is ergodic.
\end{lemma}
\begin{proof}
Given any fixed $\varphi\in C(X)$, the Birkhoff Ergodic Theorem
implies that there exist $\varphi^{*}\in L^{1}(\mu)$ and $Y\in \mathcal{B}(X)$
with $\mu(Y)=1$ such that, for any $y\in Y$,
$\lim_{n\rightarrow+\infty}(1/ n) \sum_{i=0}^{n-1}\varphi \circ T^{i}(y)= \varphi^{*}(y)$.
Since $(T, \mu)$ satisfies LDT-III', for any $N\in \mathbb{N}$ there
exists an increasing sequence $\{L_{k}^{(N)}\}_{k=1}^{+\infty}$ such that,
for any $k\in \mathbb{N}$,
$$
\mu\left(\Sigma(L_{k}^{(N)}, 1/N)=\left\{x\in X: \left|\frac{1}{L_{k}^{(N)}}\sum_{i=0}^{L_{k}^{(N)}-1}
\varphi \circ T^{i}(x)-\int_{X} \varphi \mathrm{d}\mu\right|>\frac{1}{ N}\right\}\right)\leq \frac{1}{2^{k}}.
$$
Set $\Sigma=\bigcap_{N=1}^{+\infty}\bigcup_{m=1}^{+\infty}\bigcap_{k=m}^{+\infty}X\setminus \Sigma(L_{k}^{(N)}, 1/N)
=X\setminus \bigcup_{N=1}^{+\infty}\bigcap_{m=1}^{+\infty}\bigcup_{k=m}^{+\infty} \Sigma(L_{k}^{(N)}, 1/N)$,
and take $\Omega=\Sigma\cap Y$.
Then, $\mu(\Omega)=1$. According to the construction of $\Omega$, it follows that,
for any $x\in \Omega$ and any $N\in \mathbb{N}$, there exists $m\in \mathbb{N}$ such
that for any $k\geq m$,
$$
\left|\frac{1}{L_{k}^{(N)}}\sum_{i=0}^{L_{k}^{(N)}-1}
\varphi \circ T^{i}(x)-\int_{X} \varphi \mathrm{d}\mu\right|\leq\frac{1}{N}.
$$
Combining this with $\lim_{k\rightarrow +\infty}(1/ L_{k}^{(N)})\sum_{i=0}^{L_{k}^{(N)}-1}\varphi \circ T^{i}(y)
=\varphi^{*}(x)$
yields that, for any $x\in \Omega$,
$$
\lim_{n\rightarrow+\infty}\frac{1}{n}\sum_{i=0}^{n-1}\varphi \circ T^{i}(x)
=\varphi^{*}(x)=\int_{X} \varphi \mathrm{d}\mu.
$$
The proof is then completed by Lemma \ref{ergodic-char}.
\end{proof}

\begin{theorem}\label{LDT-III'}
Let $(X, T)$ be a dynamical system and $\mu\in M(X, T)$. Then, the
following statements are equivalent:
\begin{enumerate}
  \item $(T, \mu)$ satisfies LDT-II'.
  \item $(T, \mu)$ satisfies LDT-III'.
  \item $\mu$ is ergodic.
\end{enumerate}
\end{theorem}
\begin{proof}
They follow from Theorem \ref{LDT-II'} and Lemma \ref{Lemma 6.1}.
\end{proof}


\begin{corollary}\label{Corollary 3-3}
Let $(X, T)$ be a uniquely ergodic dynamical system and $\mu\in M(X)$. Then,
the following statements are equivalent:
\begin{enumerate}
  \item $(T, \mu)$ satisfies LDT.
  \item $(T, \mu)$ satisfies LDT-I.
  \item $(T, \mu)$ satisfies LDT-II.
  \item $T$ is minimal and $\mu\in M(X, T)$.
  \item $(T, \mu)$ satisfies LDT-III and $\mu\in M(X, T)$.
\end{enumerate}
\end{corollary}

\begin{remark}
Applying Theorem \ref{Theorem 1}, it is not difficult to check
that Theorem \ref{sen-Th-1}, Corollary \ref{sen-Th}, Corollary \ref{Li},
and Corollary \ref{Corollary 3.3} all hold for LDT-III.
\end{remark}

Although Theorem \ref{LDT-III'} proves that for a measure-preserving transformation,
LDT-II' is equivalent to LDT-III', we do not know if this holds for general
transformations, because we can not answer the question as if the probability measure $\mu$ is
$T$-invariant when $(T, \mu)$ satisfies LDT-III'. So, the following question is posed.

\begin{question}
  Is there a transformation of a probability space $(X, \mathcal{B}(X), \mu)$, satisfying LDT-III' (resp., LDT-II'),
  which does not satisfy LDT-II' (resp., LDT-I')?
\end{question}

\section{A remark on the central limit theorem}

This section shows
that, in contrast with the case of the large deviations theorem (see Corollary \ref{Corollary 3-3}
and Remark \ref{Remark-2}), a dynamical system satisfying the Central Limit Theorem
in the sense of Niu and Wang \cite{Niu-Wang2010} does not exist.
\begin{theorem}\label{CLT}
There exists no dynamical system satisfying the Central Limit Theorem
in the sense of Niu and Wang.
\end{theorem}
\begin{proof}
For any given dynamical system $(X, T)$, define
an observable $\varphi: X\rightarrow \mathbb{C}$
as $\varphi(x)\equiv 0$. Take a probability measure $\mu\in M(X)$.
For every interval $A\subset \mathbb{R}$, one has
\begin{eqnarray*}
&&\quad \lim_{n\rightarrow\infty}\mu \left(\left\{x\in X:
\frac{1}{\sqrt{n}}\sum_{i=0}^{n-1}\left(\varphi\circ T^{i}(x)- \int_{X}
\varphi \mathrm{d}\mu\right)\in A\right\}\right)\\
&&=\mu \left(\left\{x\in X:
0\in A\right\}\right)=
\left\{\begin{array}{cc}
1, & {\rm } \ 0\in A, \\\
0, & {\rm } \ 0\notin A,
\end{array}
\right.
=\mathbf{1}_{A}(0),
\end{eqnarray*}
where $\mathbf{1}_{A}(\cdot)$ is the characteristic function of the set $A$.
Clearly,
$$\mathbf{1}_{A}(0)=\frac{1}{\sigma
\sqrt{2\pi}}\int_{A}e^{-t^{2}/ (2\sigma^{2})}\mathrm{d}t
$$
does not hold for any interval $A\subset \mathbb{R}$. This implies that
$\varphi$ does not satisfy the Central Limit Theorem in the sense of
Niu and Wang. So, $(T, \mu)$ does not satisfy the Central Limit
Theorem in the sense of Niu and Wang.
\end{proof}
\begin{remark}
\cite[Question 4.9]{WC2014} asks if there exists a dynamical system satisfying the Central Limit Theorem?
Theorem \ref{CLT} above shows that the answer is negative.
\end{remark}


\bibliographystyle{amsplain}

\end{document}